\documentclass[12pt]{article}
\usepackage{amssymb}

\newtheorem{definition}{Definition}[section]
\newtheorem{theorem}[definition]{Theorem}
\newtheorem{lemma}[definition]{Lemma}
\newtheorem{corollary}[definition]{Corollary}
\newtheorem{remark}[definition]{Remark}
\newtheorem{notation}[definition]{Notation}

\newenvironment{proof}
{\noindent{\em Proof:\ }}{\hfill$\Box$\\}
\newenvironment{proofof}[1]
{\medskip\noindent{\em Proof of {#1}:\ }}{\hfill$\Box$\\}

\newcommand\interval[3]{{#1}\leq{#2}\leq{#3}}

\def\R{\mathbb R}
\def\Z{\mathbb Z}
\def\Q{\mathbb Q}

\newcommand\be{\begin{equation}}
\newcommand\ee{\end{equation}}

\newcommand\inv{^{-1}}
\newcommand\dist{\partial}

\begin{document}
\title{ \bf  Almost-bipartite distance-regular graphs\\
 with the $Q$-polynomial property
}
\author{Michael S. Lang and  Paul M. Terwilliger}
\date{31 August 2004}
\maketitle
\begin{abstract}
Let $\Gamma$ denote a $Q$-polynomial
distance-regular graph with diameter $D\geq4$.
Assume that the intersection numbers of $\Gamma$ satisfy $a_i=0$ for
$\interval0i{D-1}$ and $a_D\neq0$.
We show that $\Gamma$ is a polygon, a folded cube, or an Odd graph.
\end{abstract}

%*****************************************************************************
\section{Introduction}
%*****************************************************************************

\noindent
In this article we prove the following theorem.

\begin{theorem}
\label{thm:main}
Let $\Gamma$ denote a distance-regular graph with
diameter $D\geq 3$.
Assume that the intersection numbers of $\Gamma$ satisfy
$a_i=0$ for $0 \leq i\leq D-1$ and $a_D\not=0$.
Then $\Gamma$ is $Q$-polynomial if and only if 
at least one of (\ref{i:cycle})--(\ref{i:d3}) holds below.
\begin{enumerate}
\item\label{i:cycle}
$\Gamma$ is the $(2D+1)$-gon. 
\item\label{i:cubequot}
$\Gamma$ is the folded $(2D+1)$-cube.
\item\label{i:odd}
$\Gamma$ is the Odd graph on a set of cardinality $2D+1$.
\item\label{i:d3}
$D=3$ and there exist complex scalars $\beta$ and $\mu$
such that the intersection numbers of $\Gamma$ satisfy
\begin{eqnarray}
k&=&1+(\beta^2-1)\bigg(\beta(\beta+2)-(\beta+1)\mu\bigg),\label{e:k}\\
c_2&=&\mu,\\
c_3&=&-(\beta+1)(\beta^2+\beta-1-(\beta+1)\mu).\label{e:c3}
\end{eqnarray}
\end{enumerate}
\end{theorem}

\noindent
The following remarks refer to Theorem \ref{thm:main}.

\begin{remark}\label{rem:evals}
\rm
Suppose that (\ref{i:d3}) holds.
Then $\theta_0,\theta_1,\theta_2,\theta_3$ is a $Q$-polynomial ordering of
the eigenvalues of $\Gamma$, where
\begin{eqnarray}
\theta_0&=&1+(\beta^2-1)\bigg(\beta(\beta+2)-(\beta+1)\mu\bigg),\label{e:t0}\\
\theta_1&=&(\beta+1)(\beta^2+\beta-1-\beta\mu)\label{e:t1},\\
\theta_2&=&\beta^2+\beta-1-(\beta+1)\mu,\\
\theta_3&=&1-\beta-\beta^2.\label{e:t3}
\end{eqnarray}
\end{remark}

\begin{remark}\label{rem:specials}
\rm
$\Gamma$ is the 7-gon if and only if (\ref{i:d3}) holds
with $\mu=1$ and
$\beta\in\{\omega+\omega^{-1},\omega^2+\omega^{-2},\omega^3+\omega^{-3}\}$,
where $\omega$ is a primitive 7th root of unity.
$\Gamma$ is the folded 7-cube if and only if (\ref{i:d3}) holds
with $\mu=2$ and $\beta\in\{-2,2\}$.
$\Gamma$ is the Odd graph on a set of cardinality 7
if and only if (\ref{i:d3}) holds
with $\mu=1$ and $\beta=-2$.
\end{remark}

\begin{remark}\label{rem:beta}
\rm
Suppose that (\ref{i:d3}) holds but none of (\ref{i:cycle})--(\ref{i:odd}) do.
Then $\beta$ is unique, integral and less than $-2$.
We know of no graph for which this occurs.
\end{remark}

%*****************************************************************************
\section{Preliminaries}
%*****************************************************************************

Let $\Gamma=(X,R)$ denote a finite, undirected, connected graph,
without loops or multiple edges, with vertex set $X$, edge set $R$,
path-length distance function $\dist$,
and diameter $D:=\hbox{max}\{\dist(x,y):x,y\in X\}$.
Let $k$ denote a nonnegative integer.
We say $\Gamma$ is {\em regular} with {\em valency} $k$
whenever for all $x\in X$, $|\{z\in X:\dist(x,z)=1\}|=k$.
We say $\Gamma$ is {\em distance-regular} whenever
for all integers $h,i,j$ $(\interval0{h,i,j}D)$ and all $x,y\in X$ with
$\dist(x,y)=h$, the scalar
$p^h_{ij}:=|\{z\in X:\dist(x,z)=i,\dist(y,z)=j\}|$
is independent of $x$ and $y$.
For notational convenience, set
$c_i:=p^i_{1i-1}$ $(\interval1iD)$,
$a_i:=p^i_{1i}$ $(\interval0iD)$,
$b_i:=p^i_{1i+1}$ $(\interval0i{D-1})$,
$c_0:=0$ and $b_D:=0$.
For the rest of this section, suppose that $\Gamma$ is distance-regular.
We observe that $\Gamma$ is regular with valency $k=b_0$.
Further, we observe $c_i+a_i+b_i=k$ for $\interval0iD$.

We recall the Bose-Mesner algebra.
Let $\R$ denote the field of real numbers.
By $\mathit{Mat}_X(\R)$ we mean the $\R$-algebra consisting of all matrices
whose entries are in $\R$ and
whose rows and columns are indexed by $X$.
For each integer $i$ $(\interval0iD)$,
let $A_i$ denote the matrix in $\mathit{Mat}_X(\R)$ with $x,y$
entry
$$
(A_i)_{xy}=\left\{
\begin{array}{ll}
1&\hbox{if }\dist(x,y)=i,\\
0&\hbox{otherwise}
\end{array}\right.\qquad\qquad(x,y\in X).
$$
Note that $A_0=I$, the identity matrix.
Abbreviate $A:=A_1$.
We call $A$ the {\em adjacency matrix} of $\Gamma$.
Let $M$ denote the subalgebra of $\mathit{Mat}_X(\R)$ generated by $A$.
By \cite[Theorem 20.7]{biggs}, $A_0,A_1,\ldots,A_D$ is a basis for $M$.
We call $M$ the {\em Bose-Mesner algebra} of $\Gamma$.

By \cite[Theorem 2.6.1]{bcn},
$M$ has a
second
basis $E_0,E_1,\ldots,E_D$ such that
$E_iE_j=\delta_{ij}E_i$
$(\interval0{i,j}D)$.
We call $E_0,E_1,\ldots,E_D$ the {\em primitive idempotents} of $\Gamma$.
Observe that
there exists a sequence of scalars $\theta_0,\theta_1,\ldots,\theta_D$
taken from $\R$ such that
$
A=\sum_{i=0}^D\theta_iE_i.
$
We call $\theta_i$ the {\em eigenvalue} of $\Gamma$
associated with $E_i$ $(\interval0iD)$.
Note that $\theta_0,\theta_1,\ldots,\theta_D$
are distinct since $A$ generates $M$.

We recall the $Q$-polynomial property.
Let $\theta_0,\theta_1,...,\theta_D$ denote an ordering of the eigenvalues
of $\Gamma$.
We say this ordering is {\em $Q$-polynomial}
whenever there exists a sequence of real scalars
$\sigma_0,\sigma_1,...,\sigma_D$ and a sequence of polynomials
$q_0,q_1,...,q_D$ with real coefficients such that
$q_j$ has degree $j$ and
$E_j=\sum_{i=0}^Dq_j(\sigma_i)A_i$
for $\interval0jD$.
In this case, $\theta_0=k$ \cite[Theorem 8.1.1]{bcn};
we call $\theta_1$ a {\em $Q$-polynomial} eigenvalue of $\Gamma$.
We say that $\Gamma$ is {\em $Q$-polynomial} whenever there exists a
$Q$-polynomial ordering of its eigenvalues.

We recall what it means for $\Gamma$ to be bipartite or almost-bipartite.
We say $\Gamma$ is {\em bipartite} whenever $a_i=0$ for $\interval0iD$.
We say $\Gamma$ is {\em almost-bipartite} whenever
$a_i=0$ for $\interval0i{D-1}$ but $a_D\neq0$.
(In the literature, such a $\Gamma$ is also called a {\em generalized Odd
  graph} or a {\em regular thin near $(2D+1)$-gon}.)
For the rest of this section, assume that $\Gamma$ is almost-bipartite.

We recall the {\em bipartite double} $2.\Gamma$.
This graph has vertex set $\{x^+:x\in X\}\cup\{y^-:y\in X\}$.
For $x,y\in X$ and $\gamma,\delta\in\{+,-\}$,
vertices $x^\gamma$ and $y^\delta$ are adjacent in $2.\Gamma$ whenever
$x$ and $y$ are adjacent in $\Gamma$ and $\gamma\neq\delta$.
The graph $2.{\Gamma}$ is bipartite and distance-regular
with diameter $2D+1$.
Moreover, $2.\Gamma$ is an antipodal 2-cover of $\Gamma$
\cite[Theorem 1.11.1(i),(vi)]{bcn}.
The intersection numbers $k$ and $c_2$ are the same in $2.\Gamma$ as in
$\Gamma$ \cite[Proposition 4.2.2(ii)]{bcn}.
The set of eigenvalues for $2.\Gamma$ consists of the eigenvalues of $\Gamma$
together with their opposites \cite[Theorem 1.11.1(v)]{bcn}.
The concept of an AO eigenvalue was introduced in \cite{langBipIneq}.
A scalar $\theta$ is an AO eigenvalue of $2.\Gamma$ if and only if $\theta$
is a $Q$-polynomial eigenvalue of $\Gamma$ \cite[Theorem 10.4]{langBipIneq}.

%*****************************************************************************
\section{Setup}
%*****************************************************************************

Our strategy for proving Theorem \ref{thm:main} is to assume that
(\ref{i:cycle})--(\ref{i:odd}) do not hold and then prove that
(\ref{i:d3}) must.

\begin{lemma}\label{l:1qpoly}
Let $\Gamma$ denote an almost-bipartite distance-regular graph with
diameter $D\geq3$.
Suppose that $\Gamma$ is $Q$-polynomial but not as in Theorem
\ref{thm:main}(\ref{i:cycle})--(\ref{i:odd}).
Then $\Gamma$ has a unique $Q$-polynomial eigenvalue.
\end{lemma}

\begin{proof}
Since $\Gamma$ is not a polygon, it has valency $k\geq3$.
Observe that $2.\Gamma$ has diameter at least 7.
Suppose that $\Gamma$ has at least two $Q$-polynomial eigenvalues.
Then $2.\Gamma$ has at least two AO eigenvalues.
Applying \cite[Theorem 16.2]{langBipIneq}, we find that $2.\Gamma$
is the $(2D+1)$-cube.
Thus $\Gamma$ is the folded $(2D+1)$-cube, contradicting the
assumption that Theorem \ref{thm:main}(\ref{i:cubequot}) does not hold.
\end{proof}

\noindent
For the rest of this article, we use the following notation.

\begin{notation}\label{n:gamma}
Let $\Gamma=(X,R)$ denote an almost-bipartite distance-regular graph
with diameter $D\geq 3$.
Assume that $\Gamma$ is $Q$-polynomial but not as in
Theorem \ref{thm:main}(\ref{i:cycle})--(\ref{i:odd}).
Let $\theta_0,\theta_1,\ldots,\theta_D$ denote the eigenvalues of
$\Gamma$ in their $Q$-polynomial order.
Set
\be\label{e:beta}
\beta:=\frac{\theta_0-\theta_3}{\theta_1-\theta_2}-1.
\ee
\end{notation}

%*****************************************************************************
\section{Parameters}
%*****************************************************************************

In this section, we recall some formulae for the intersection numbers and
eigenvalues of $\Gamma$.

\begin{lemma}\label{lem:parms}
\cite[Lemma 15.2, Corollaries 15.4, 15.7, 15.8, Theorem 15.5]{cmtAlmBip}
With reference to Notation \ref{n:gamma},
there exist complex scalars $q$ and $s$ such that
the intersection numbers and eigenvalues of $\Gamma$ satisfy
\begin{eqnarray}
k&=&h(1+sq),\label{eq:k}\\
c_i&=&\frac{h(1-q^i)(1+sq^{2D+2-i})}{q^i(q^{2D-2i+1}-1)}
\qquad(1\leq i\leq D)\label{eq:ci},\\
&&\nonumber\\
\theta_i&=&hq^{-i}(1+sq^{2i+1})\qquad(\interval0iD),\label{eq:thform}
\end{eqnarray}
where
\be
h = \frac{q-q^{2D}}{(q-1)(1+sq^{2D+1})}.\label{eq:h}
\ee
Moreover,
\begin{eqnarray}
q&\neq&0,\label{e:hq}\\
q^i&\neq&1\qquad(\interval1i{2D}),\label{eq:qrest}\\
sq^i&\neq&1\qquad(\interval2i{2D}),\\
sq^i&\neq&-1\qquad(\interval1i{2D+1})\label{eq:smore}.
\end{eqnarray}
\end{lemma}

\begin{corollary}
With reference to Lemma \ref{lem:parms}, we have
\be\label{e:betaq}
\beta=q+q\inv
\ee
and
\be\label{eq:thd}
\theta_D=\frac{q^{1-D}-q^D}{q-1}.
\ee
\end{corollary}

%*****************************************************************************
\section{Restrictions}
%*****************************************************************************

Throughout this section we refer to Notation \ref{n:gamma} and Lemma
\ref{lem:parms}.
We obtain restrictions on the parameters $q$ and $s$.
Let $\Z$ denote the ring of integers;
let $\Q$ denote the field of rational numbers.

\begin{lemma} 
\label{lem:eigint}
We have $\theta_i\in\Z$ for $0\leq i\leq D$.
\end{lemma}

\begin{proof}
Suppose that there exists an integer $i$ $(0 \leq i \leq D)$ such
that $\theta_i\not\in\Z$.
Then $\Gamma$ has a second $Q$-polynomial eigenvalue by
\cite[p. 360]{bannai}.
This contradicts Lemma \ref{l:1qpoly}.
\end{proof}

\begin{lemma} 
\label{lem:betaz}
We have $q^i+q^{-i} \in\Z$ for each positive integer $i$.
In particular, $\beta \in \Z$.
\end{lemma}

\begin{proof}
Define polynomials $T_0, T_1, ...$ in a variable $x$ by
$T_0=2$,
$T_1=x$,
$T_{i+1}=xT_i-T_{i-1}$ $(i\geq1)$.
We routinely find that for $i\geq1$,
\be\label{e:tiobs}
T_i\in\Z[x],\quad
T_i\hbox{ is monic},\quad
\hbox{and}\quad
q^i+q^{-i}=T_i(\beta).
\ee

To finish the proof it suffices to show
$\beta \in \Z$.
To do this we show $\beta \in \Q$ and $\beta$ is an algebraic integer.

By (\ref{e:beta}) and Lemma \ref{lem:eigint} we find $\beta\in\Q$. 

We now show that $\beta$ is an algebraic integer.
The right-hand side of
(\ref{eq:thd}) is equal to $-\sum_{i=1-D}^{D-1} q^i$.
By this and (\ref{e:tiobs}),
we find that $\beta$ is a root of a monic polynomial with coefficients in $\Z$.
It follows that $\beta$ is an algebraic integer.

We have now shown $\beta \in\Q$ and $\beta$ is an algebraic integer,
so $\beta \in \Z$.
The result follows.
\end{proof}

\begin{lemma}
\label{lem:qreal}
We have $|\beta|>2$.
Moreover, $q \in \R$.
\end{lemma}

\begin{proof}
Suppose $|\beta|\leq2$. 
Since $\beta \in \Z$ by Lemma \ref{lem:betaz},
we find that $|\beta|$ is 0, 1 or 2.
We now use (\ref{e:betaq}).
If $|\beta|=0$ then $q^4=1$.
If $|\beta|=1$ then $q^6=1$.
If $|\beta|=2$ then $q^2=1$.
Each of these contradicts (\ref{eq:qrest}), so the result follows.
\end{proof}

\begin{lemma}\label{l:q2>1}
We may assume $q^2>1$.
\end{lemma}

\begin{proof}
By (\ref{e:hq}), we have $q^2\neq0$.
By (\ref{eq:qrest}), we have $q^2\neq1$.
We now consider two cases.

First suppose $s=0$.
If $-1<q<0$ then using (\ref{eq:ci}) we find $c_2<0$.
If $0<q<1$ then using (\ref{eq:k}) we find $k<0$.
Each of these is a contradiction, so $q^2>1$ as desired.

Now suppose $s\neq0$.
If $q^2<1$, replace $q$ by $q\inv$ and $s$ by $s\inv$.
In light of (\ref{eq:h}), these substitutions leave
(\ref{eq:k})--(\ref{eq:thform}) unchanged.
Moreover, $q^2>1$ as desired.
\end{proof}

\noindent
Consider the quantity
\be\label{e:eta}
\eta := -\frac{(q^2+1)(q^{2D}-q^3)}{q^{2D}-q^5}.
\ee
We show $\eta$ to be an integer.
To do this, we use the fact $s^2q^{2D+3}\neq1$.
We obtain this fact using the following two lemmas. 

\begin{lemma}\label{l:curtinIneq}
For $1 \leq i \leq D$ we  have 
$(c_2-1)\theta_i^2\neq(k-c_2)(k-2)$.
\end{lemma}

\begin{proof}
Suppose that there exists an integer $i$ $(\interval1iD)$ such that
$(c_2-1)\theta_i^2=(k-c_2)(k-2)$.
We mentioned earlier that $\theta_i$ is an eigenvalue of $2.\Gamma$ and
that the intersection numbers $k$ and $c_2$ are the same in $2.\Gamma$ as
in $\Gamma$.
Now by \cite[Theorem 25]{curtin2homo}, we find that $2.\Gamma$
is 2-homogeneous in the sense of Nomura \cite{nomuraHomo}.
By assumption, $2.\Gamma$ is not a cube.
Now by \cite[Theorem 1.2]{nomuraSpin}, the diameter of $2.\Gamma$ is at most 5.
Since this diameter is $2D+1$, we find $D\leq 2$, which is a contradiction.
\end{proof}

\begin{lemma}\label{l:curtinEval}
We have
$$(c_2-1)\theta_D^2-(k-c_2)(k-2)=
\frac{(q^{2D}-1)(q^{2D}-q^2)(q^{2D}-q)^2(s^2q^{2D+3}-1)}
{q^{2D}(q-1)^2(q^{2D}-q^3)(1+sq^{2D+1})^2}.
$$
\end{lemma}
\begin{proof}
Use Lemma \ref{lem:parms}.
\end{proof}

\begin{corollary}\label{cor:s2q2d+3}
We have
$s^2q^{2D+3}\neq1$.
\end{corollary}
\begin{proof}
Combine Lemmas \ref{l:curtinIneq} and \ref{l:curtinEval}.
\end{proof}

\noindent
Before proceeding, we recall the local graph $\Gamma^2_2$.

\begin{definition}\rm\label{def:g22}
Fix a vertex $x\in X$.
The corresponding {\em local graph} $\Gamma^2_2$ is
the graph with vertex set $\{y\in X:\dist(x,y)=2\}$,
where vertices $y$ and $z$ are adjacent in $\Gamma^2_2$ whenever $\dist(y,z)=2$
in $\Gamma$.
\end{definition}

\begin{lemma}
\label{lem:loc}
Fix a vertex $x\in X$ and
let $\Gamma^2_2$ denote the corresponding local graph
from Definition \ref{def:g22}.
Then the scalar $\eta$ from (\ref{e:eta}) is an eigenvalue of $\Gamma^2_2$.
Moreover, $\eta$ is an algebraic integer.
\end{lemma}

\begin{proof}
Our argument uses the subconstituent algebra of $\Gamma$.
This object is introduced in \cite{terwSubconstituentI}.
We refer the reader to that paper and its continuations
\cite{terwSubconstituentII} and \cite{terwSubconstituentIII}
for background and definitions.

Let $T=T(x)$ denote the subconstituent algebra of $\Gamma$ with respect to $x$.
By \cite[Theorem 14.1]{cmtAlmBip} we find that, up to isomorphism,
there exists at most one
irreducible $T$-module with endpoint 2, dual endpoint 2 and diameter $D-2$.
By \cite[Example 16.9(iv)]{cmtAlmBip}
the multiplicity with which this module appears in the standard module is
$$
\frac{(q^{2D}-1)(q^{2D}-q^2)(1+sq)(1+sq^4)(s^2q^{2D+3}-1)}
{q(q+1)(q-1)^2(s^2q^{2D+4}-1)(1+sq^{2D})(1+sq^{2D+1})}.
$$
This number is nonzero by (\ref{e:hq}), (\ref{eq:qrest}), (\ref{eq:smore}) and Corollary
\ref{cor:s2q2d+3}.
Therefore, this module exists.

Let $W$ denote an
irreducible $T$-module with endpoint 2, dual endpoint 2 and diameter $D-2$.
The dimension of $E^*_iW$ is 1 for $\interval2iD$
\cite[(18), Lemma 10.3]{cmtAlmBip}.
By construction, $E^*_2W$ is an eigenspace for $E^*_2A_2E^*_2$.
By \cite[Definition 8.2]{cmtAlmBip} and using $A_2=(A^2-kI)/c_2$,
we find that the corresponding eigenvalue is
\be\label{e:eval}
\frac{c_1(W)b_0(W)-k}{c_2},
\ee
where $c_1(W)$ and $b_0(W)$ are intersection numbers of $W$.
Evaluating (\ref{e:eval}) using \cite[Theorem 15.5]{cmtAlmBip},
we find that it is equal to $\eta$.
We conclude that $\eta$ is an eigenvalue of $\Gamma^2_2$.

Now $\eta$ is an algebraic integer since it is an eigenvalue of a graph.
\end{proof}

\noindent
We show $\eta \in \Z$.
To do this we use the following result.

\begin{lemma}
We have 
\begin{eqnarray}
\label{eq:exp}
\eta +\beta^2-1 =\frac{q^{2D}-q^9}{q^{2D+2}-q^7}.
\end{eqnarray}
\end{lemma}
\begin{proof}
Use (\ref{e:betaq}) and (\ref{e:eta}).
\end{proof}

\begin{lemma}
\label{lem:etaint}
We have $\eta \in \Z$.
\end{lemma}

\begin{proof}
First assume $D=3$.
Then $\eta=-\beta(\beta+1)$ by (\ref{e:betaq}) and (\ref{e:eta}).
Thus $\eta\in\Z$ by Lemma \ref{lem:betaz}.
Now assume $D\geq 4$.
Observe by Lemma
\ref{lem:loc} that $\eta$ is an algebraic integer.
We show
 $\eta\in\Q$.
Observe that the right-hand side of (\ref{eq:exp}) is equal to
$-(\beta+1)^{-1}$ for $D=4$ and
\[\frac{\sum_{i=5-D}^{D-5}q^i}{\sum_{i=3-D}^{D-3}q^i}\]
for $D\geq 5$.
By this and
Lemma
\ref{lem:betaz}
we find that the right-hand side of
(\ref{eq:exp}) is in $\Q$.
By this and since $\beta \in \Z$ we find
$\eta \in \Q$. 
Now  $\eta \in \Q$ and $\eta$ is an algebraic integer
so 
$\eta \in \Z$.
\end{proof}

%*****************************************************************************
\section{Proof}
%*****************************************************************************

In this section we prove Theorem \ref{thm:main} and the associated remarks.

\begin{proofof}{Theorem \ref{thm:main}}
Assume that $\Gamma$ is $Q$-polynomial but
none of (\ref{i:cycle})--(\ref{i:odd}) hold.
We show that $\Gamma$ satisfies (\ref{i:d3}).

We first show $D=3$.
On the contrary, suppose $D\geq 4$.
For notational convenience, abbreviate $\xi:=\eta+\beta^2-1$.
Recall
$\beta\in\Z$ by Lemma \ref{lem:betaz} and
$\eta\in\Z$ by Lemma \ref{lem:etaint}, so $\xi\in\Z$.
By (\ref{e:hq}), (\ref{eq:qrest}) and (\ref{eq:exp})
we find
$\xi\neq0$.
Thus $|\xi|\geq1$.
Evaluating $\xi^2-1$ using
(\ref{eq:exp}) and simplifying, we find
$(q^4-1)(q^{14}-q^{4D})\geq0$.
But since $q^2>1$ by Lemma \ref{l:q2>1},
we find
$(q^4-1)(q^{14}-q^{4D})<0$,
for a contradiction.
We have now shown $D=3$.

Evaluating (\ref{eq:k})--(\ref{eq:h}) using $D=3$ and
$\beta=q+q\inv$, we routinely obtain (\ref{e:k})--(\ref{e:c3})
and (\ref{e:t0})--(\ref{e:t3}).

We have proved the theorem in one direction.
We now show the converse.
First assume that $\Gamma$ satisfies one of (\ref{i:cycle})--(\ref{i:odd}).
That $\Gamma$ is $Q$-polynomial is well known
\cite[Corollary 8.5.3]{bcn}.
Now assume that $\Gamma$ satisfies (\ref{i:d3}).
We routinely find that the eigenvalues of $\Gamma$,
in a $Q$-polynomial order, are given by (\ref{e:t0})--(\ref{e:t3}).
\end{proofof}

\noindent
Remarks \ref{rem:evals} and \ref{rem:specials} are verified routinely.

\begin{proofof}{Remark \ref{rem:beta}}
The $Q$-polynomial ordering of the eigenvalues
is unique by Lemma \ref{l:1qpoly},
so $\beta$ is unique by (\ref{e:beta}).
It is an integer by Lemma \ref{lem:betaz}.

We show $\beta<-2$.
Recall that $\theta_1$ is a $Q$-polynomial eigenvalue of $\Gamma$ and thus
is an AO eigenvalue of the bipartite double $2.\Gamma$.
Observe that $2.\Gamma$ has diameter 7.
Since $2.\Gamma$ is bipartite, we see by \cite[p. 82]{bcn} that half of the
eigenvalues of $2.\Gamma$ are positive and half are negative.
By \cite[Lemma 13.5]{langBipIneq}, we find that $\theta_1$ is the fifth- or
seventh-largest of the eight eigenvalues of $2.\Gamma$.
Thus $\theta_1<0$.

Recall $|\beta|>2$ by Lemma \ref{lem:qreal};
thus $\beta^2+\beta-1>0$.
Observe
$b_2=(\beta^2+\beta-1)(\beta^2+\beta-1-\beta\mu)$
by (\ref{e:k}).
By this and since $b_2>0$, we find $\beta^2+\beta-1-\beta\mu>0$.
Combining this with (\ref{e:t1}), we find $\beta+1<0$.
In particular, $\beta<0$.
Now $\beta<-2$ in view of Lemma \ref{lem:qreal}.
\end{proofof}

\noindent
{\bf Acknowledgment.}
The first author was partially supported by a Bradley University Research
Excellence Committee award.

%*****************************************************************************
% Bibliography
%*****************************************************************************

\end{document}